\documentclass[12pt]{amsart}
\usepackage{graphicx}
\usepackage{color}
\usepackage[latin1]{inputenc}
\usepackage{amsmath}
\usepackage{amssymb}
\newtheorem{teo}{Theorem}[section]
\newtheorem{cor}[teo]{Corollary}
\newtheorem{lema}[teo]{Lemma}
\newtheorem{prop}[teo]{Proposition}
\theoremstyle{definition}

\theoremstyle{remark}
\newtheorem{rem}[teo]{Remark}

\numberwithin{equation}{section}

\newcommand{\To}{\longrightarrow}

\newcommand{\N}{\mathbb{N}}
\newcommand{\R}{\mathcal{R}}
\newcommand{\I}{\mathcal{I}}
\newcommand{\C}{\mathbb{C}}
\newcommand{\D}{\mathbf{D}}
\newcommand{\B}{\mathcal{B}}
\newcommand{\dem}{\noindent{\textsf{Proof.} }}

\renewcommand{\phi}{\varphi}
\newcommand{\eps}{\varepsilon}
\def\ba{\begin{eqnarray}}
\def\ea{\end{eqnarray}}
\def\ban{\begin{eqnarray*}}
\def\ean{\end{eqnarray*}}

\author[Miralles]{Alejandro Miralles}\address{Alejandro Miralles.
 Departament de Matemàtiques and IMAC, Universitat Jaume I, Castelló (Spain). \emph{e}.mail: mirallea@uji.es}
 
\thanks{Supported by PGC2018-094431-B-100 (MICINN. Spain) and 8059/2019 (UJI)}


\begin{document}

\begin{abstract}
Let $B_E$ be the open unit ball of a complex finite or infinite dimensional Hilbert space. If $f$ belongs to the space $\mathcal{B}(B_E)$ of Bloch functions on $B_E$, we prove that the dilation map given by $x \mapsto (1-\|x\|^2) \mathcal{R} f(x)$ for $x \in B_E$, where $\mathcal{R} f$ denotes the radial derivative of $f$, is Lipschitz continuous with respect to the pseudohyperbolic distance $\rho_E$ in $B_E$, which extends to the finite and infinite dimensional setting the result given for the classical Bloch space $\B$. In order to provide this result, we will need to prove that $\rho_E(zx,zy) \leq |z| \rho_E(x,y)$ for $x,y \in B_E$ under some conditions on $z \in \mathbb{C}$. Lipschitz continuity of $x \mapsto (1-\|x\|^2) \mathcal{R} f(x)$ will yield some applications which also extends classical results from $\mathcal{B}$ to $\mathcal{B}(B_E)$. On the one hand, we supply results on interpolating sequences for $\mathcal{B}(B_E)$: we show that it is necessary for a sequence in $B_E$ to be separated in order to be interpolating for $\mathcal{B}(B_E)$ and we also prove that any interpolating sequence for $\mathcal{B}(B_E)$ can be slightly perturbed and it remains interpolating. On the other hand, after a deep study of the automorphisms of $B_E$, we provide necessary and suficient conditions for a composition operator on $\mathcal{B}(B_E)$ to be bounded below.
\end{abstract}


\title[Lipschitz continuity of the dilation of Bloch functions]{Lipschitz continuity of the dilation of Bloch functions on the unit ball of a Hilbert space and applications}

\maketitle 

\section{Introduction and background} \label{section1}


Along this work, $E$ will denote a finite or infinite dimensional complex Hilbert space. Its open unit ball will be denoted by $B_E$. After some background, we will study in Section \ref{section2} the boundness of:
\ban \label{exp}
\frac{\rho_E(zx,zy)}{|z|\rho_E(x,y)}
\ean 

\noindent for $z \in \C$ and $x,y \in B_E$ such that $zx,zy \in B_E$. First we will show that, in general, this expression is unbounded. Nevertheless, we will show that if $|z|$ is bounded above by:
$$\frac{1+\max\{\|x\|,\|y\|\}}{2 \max\{\|x\|,\|y\|\}},$$

\noindent then the expression above is bounded by $2$ and this bound is the best possible. We first prove the case when we deal with $E=\C$ and then with any finite or infinite dimensional Hilbert space $E$. At the end of this section, we will extend this result to the case when we deal with the Banach space $C_0(S)$.

In Section \ref{section4}, we deal with functions $f$ which belong to the Bloch space $\B(B_E)$, that is, the space of Bloch functions on the unit ball $B_E$ of $E$. First, as a consequence of the boundness of (\ref{exp}), we show that the dilation map $x \mapsto (1-\|x\|^2) \R f(x)$ for $x \in B_E$ is Lipschitz with respect to the pseudohyperbolic distance in Subsection \ref{subsec31}, extending to the finite and infinite dimensional setting the result given by Attele in \cite{At} and improved by Xiong in \cite{Xi} for the classical Bloch space $\B$. Hence, we derive some results about interpolating sequences for $\B(B_E)$ in Subsection \ref{subsec42}. Indeed, we supply a proof that these sequences are separated for the pseudohyperbolic distance. We also prove that interpolating sequences can be slightly perturbed and they remain interpolating, which also extends the result for $\B$ given in \cite{At}. Finally, in Subsection \ref{subsec43} we first will study several properties of automorphisms of $B_E$ and this will permit us to give necessary and sufficient conditions for a composition operator on $\B(B_E)$ to be bounded below, which extends the results given in the one-dimensional case in \cite{Ch} to the finite and infinite dimensional setting.

\subsection{The pseudohyperbolic and hyperbolic distance.} Let $\D$ be the open unit disk of the complex plane $\C$. Recall that the pseudohyperbolic distance for $z,w \in \D$ is given by:
$$\rho(z,w)=\left| \frac{z-w}{1-\bar{z}w} \right|.$$

Let $X$ be a complex Banach space and let $B_X$ be its open unit ball. Recall that $f: B_X \to \C$ is said to be holomorphic (or analytic) if it is Fréchet differentiable for any $x \in B_X$ (see \cite{Mu} for further information). For any $x,y \in B_X$, the pseudohyperbolic distance $\rho_X(x,y)$ is given by:
\ban
\rho_X (x,y)=\sup \{ \rho(f(x),f(y)) : f \in H^{\infty}(B_X), \|f\|_{\infty} \leq 1 \}, \ean
where $H^{\infty}(B_X)$ is the space of bounded holomorphic functions on $B_X$ which become a Banach space (a uniform algebra, indeed) endowed with the sup-norm. The hyperbolic distance for $x,y \in B_X$ is given by:
$$\beta_X(x,y)=\frac{1}{2} \log \left( \frac{1+\rho_X(x,y)}{1-\rho_X(x,y)} \right).$$

\subsection{Automorphisms and pseudohyperbolic distance on $B_E$}

If we deal with a complex Hilbert space $E$, we will denote by $Aut(B_E)$ the space of automorphisms of $B_E$, that is, the maps $\phi: B_E \to B_E$ which are bijective and bianalytic. These automorphisms are well-known (see \cite{GR}) and we will need to make use of them during the sequel. For any $x \in B_E$, the automorphism $\phi_x: B_{E} \To B_{E}$ is defined according to:
\begin{eqnarray} \label{automorph}
\phi_{x}(y)= (s_x Q_x+P_x)(m_x(y))
\end{eqnarray}

\noindent where $s_x=\sqrt{1-\|x\|^2},$  $m_x: B_{E} \To B_{E}$ is  the analytic self-map: $$m_x(y)=\frac{x-y}{1- \langle y,x \rangle},$$ \noindent $P_x: E \To E$ is the orthogonal projection along the one-dimensional subspace spanned by $x$, that is:
$$P_x(y)=\frac{\langle y, x \rangle}{\langle x,x \rangle} x$$
\noindent and $Q_x: E \To E$  is the orthogonal complement $Q_x=Id_E-P_x$, where $Id_E$ denotes the identity operator on $E$. It is clear that $\phi_x(0)=x$ and $\phi_x(x)=0$. The automorphisms of the unit ball $B_{E}$ turn to be compositions of these $\phi_x$ with  unitary transformations $U$ of $E$. \medskip

It is well-known (see \cite{GR}) that the pseudohyperbolic distance on $B_E$ is given by:
\begin{eqnarray} \label{metrica hilb}
\rho_{E}(x,y) = \|\phi_{y}(x)\| \mbox{ for any } x, y \in B_{E}.
\end{eqnarray}
\noindent and:
\ba \label{pseudoh}
\rho_E(x,y)^2=1-\frac{(1-\|x\|^2)(1-\|y\|^2)}{|1-\langle x,y \rangle|^2}.
\ea

\subsection{The Bloch space} The classical Bloch space $\B$ is the set of holomorphic functions $f: \D \to \C$ such that $\|f\|_B=\sup_{z \in \D} (1-|z|^2)|f'(z)|$ is bounded. This supremum defines a semi-norm which becomes a norm by adding up a constant: $|f(0)|+\sup_{z \in \D} (1-|z|^2)|f'(z)|$. Hence, $\B$ becomes a complex Banach space. The semi-norm $\| \cdot \|_B$ is invariant by automorphisms, that is, $\|f \circ \phi\|_B=\|f\|_B$ for any $f \in \B$ and $\phi: \D \to \D$ an automorphism of $\D$. \medskip

Timoney extended Bloch functions if we deal with a finite dimensional Hilbert space (see \cite{T80}). Blasco, Galindo and Miralles extended them to the infinite dimensional setting (see \cite{BGM}). If we deal with a complex finite or infinite dimensional Hilbert space $E$, the analytic function $f: B_E \to \C$ is said to belong to the Bloch space $\B(B_E)$ if:
$$\|f\|_\B=\sup_{x \in B_E} (1-\|x\|^2) \| \nabla f(x)\| < +\infty,$$
\noindent where it is clear that $\nabla f(x)$ is the derivative $f'(x)$ or, equivalently, if:
$$\|f\|_\R=\sup_{x \in B_E} (1-\|x\|^2) \| R f(x)\| < +\infty,$$
\noindent where $\R f(x)$ is the radial derivative of $f$ at $x$ given by $\R f(x)=\langle x, \overline{\nabla f (x)}\rangle$. These semi-norms are equivalent to the following one:
\ba \label{invar}
\|f\|_\I=\sup_{x \in B_E} \| \widetilde{\nabla} f (x) \|,
\ea
\noindent where $\widetilde{\nabla} f (x)$ denotes the invariant gradient of $f$ at $x$ which is given by $\widetilde{\nabla} f (x) =\nabla (f \circ \phi_x)(0)$, where $\phi_x$ is the automorphism given in (\ref{automorph}). \medskip

The three semi-norms $\| \cdot \|_\B, \| \cdot \|_\R$ and $\| \cdot \|_\I$ define equivalent Banach space norms-modulo the constant functions- in $\B(B_E)$ (see \cite{BGM}). In particular, there exists a constant $A_0 >0$ such that:
\ba \label{seminormas}
\|f\|_\R \leq \|f\|_\B \leq \|f\|_\I \leq A_0 \|f\|_R. 
\ea

Hence, the space $\B(B_E)$ can be endowed with any of the norms:
$$\| \cdot \|_{\B-Bloch}=|f(0)|+ \|\cdot\|_\B$$
or:
$$\| \cdot \|_{\R-Bloch}=|f(0)|+ \|\cdot\|_\R$$
or:
$$\| \cdot \|_{\I-Bloch}=|f(0)|+ \|\cdot\|_\I$$
and $\B(B_E)$ becomes a Banach space. We will make use of these three semi-norms and norms along the sequel. \noindent We will also make use of this result, which states that Bloch functions on $B_E$ are Lipschitz with respect to the hyperbolic distance (see \cite{BGLM}):
\begin{prop} \label{prop_lip}
	Let $E$ be a complex Hilbert space and let $f \in \B(B_E)$. Then for any $x,y \in B_E$:
	$$|f(x) -f(y)| \leq \|f\|_\I \beta_E(x,y).$$
\end{prop}

\section{Inequalities with the pseudohyperbolic distance} \label{section2}

Let $X$ be a complex Banach space. If $\phi: B_X \to B_X$ is an analytic self-map, it is well-known that $\rho_X(\phi(x),\phi(y)) \leq \rho_X (x,y)$ for any $x,y \in B_X$ and the equality is attained if and only if $\phi$ is an automorphism of $B_X$. Hence, if we consider $z \in \C$, $|z| \leq 1$ and $x,y \in B_X$, it is clear that:
\ban \label{condi}
\frac{\rho_X(zx,zy)}{\rho_X(x,y)} \leq 1
\ean
\noindent since the map $\phi: B_X \to B_X$ given by $\phi(x)=z x$ is analytic on $B_X$. However, this situation changes dramatically if we consider the expression:
\ba \label{condi2}
\frac{\rho_X(zx,zy)}{|z|\rho_X(x,y)}
\ea
for any $z \in \C$ such that $zx,zy \in B_X$. We show that, in general, this expression is unbounded. Anyway, if we deal with $X$ a complex Hilbert space or $C_0(S)$ and $z \in \C$ satisfies:
$$|z| \leq \frac{1+\max\{\|x\|,\|y\|\}}{2 \max\{\|x\|,\|y\|\}},$$ then expression (\ref{condi2}) is bounded by $2$ and this will permit us to provide several applications in Section \ref{section4}.

\subsection{Unboundness}

In this section we prove that expression (\ref{condi2}) is unbounded in general.

\begin{prop} \label{counterex}
Let $E$ be a complex Hilbert space. There exist a sequence $(z_n) \subset \C$, $x \in B_E$ and a sequence $(y_n) \subset B_E$ such that $z_n x ,z_n y_n \in B_E$ but:
$$\frac{\rho_E(z_n x,z_n y_n)}{|z_n| \rho_E(x,y_n)}$$
\noindent is unbounded.
\end{prop}

\dem  We prove it for $E=\C$. Take for instance $x=1/2$, $y_n=1/2-\frac{1}{n}$ and $z_n =2-\frac{1}{n}$. It is clear that $|z_n x|<1$ and $|z_n y_n| <1$. However:
\begin{eqnarray*}
\rho(z_n x,z_ny_n)=\frac{\frac{1}{n}\left(2-\frac{1}{n} \right)}{1-\frac{1}{2}\left(\frac{1}{2}-\frac{1}{n}\right)\left(2-\frac{1}{n}\right)^2}=\frac{\frac{2n-1}{n^2}}{\frac{2-9n+12n^2}{4n^3}}=\frac{4n(2n-1)}{12n^2-9n+2}
\end{eqnarray*}
\noindent and:
\begin{eqnarray*}
|z_n|\rho( x,y_n)=\left( 2-\frac{1}{n}\right)\frac{\frac{1}{n}}{1-\frac{1}{2} \left(\frac{1}{2}-\frac{1}{n}\right)}=\frac{\frac{2n-1}{n^2}}{\frac{3n+2}{4n}}=\frac{4(2n-1)}{3n^2+2n}.
\end{eqnarray*}
\noindent Hence:
\begin{eqnarray*}
\frac{\rho(z_nx,z_ny_n)}{|z_n|\rho(x,y_n)}=\frac{\frac{4n(2n-1)}{12n^2-9n+2}}{\frac{4(2n-1)}{3n^2+2n}}=\frac{3n^3+2n^2}{12n^2-9n+2}
\end{eqnarray*}

\noindent which is clearly unbounded since it tends to $\infty$ when $n \rightarrow \infty$. The result remains true if we deal with any complex Hilbert space $E$ since we can take $x_0 \in E$ such that $\|x_0\|=1$ and take $u=\frac{1}{2} x_0$, $v_n= y_n x_0$. We have that:
\begin{eqnarray*}
\rho_E(z_n u,z_n v_n)^2=1-\frac{(1-\|u_n\|^2)(1-\|v_n\|^2)}{|1-\langle u_n,v_n \rangle|^2}= \\ 1-\frac{(1-\|z_n x\|^2)(1-\|z_n y_n\|^2)}{|1-\langle z_n x,z_n y_n \rangle|^2}=\rho_E(z_n x, z_n y_n)^2
\end{eqnarray*}
\noindent and similarly we have $\rho_E(u,v_n)=\rho_E(x,y_n)$ so we apply the case $E=\C$ and we are done. \qed \bigskip

An easy consequence is a well-known result: the pseudohyperbolic distance cannot be extended to a norm on $E$ since:
$$\frac{\rho_E(zx,zy)}{|z|\rho_E(x,y)}$$
\noindent is unbounded, so $\rho_E(zx,zy) \neq |z| \rho_E(x,y)$. \bigskip

\subsection{Boundness} 

The main result of this section is Theorem \ref{teo} which states that under condition (\ref{condi}), then the expression (\ref{exp}) is bounded and the best bound possible is given by $2$. The following lemma will be used to prove this result:
\begin{lema} \label{lema1}
Let $E$ be a finite or infinite dimensional Hilbert space, $z \in \C$ and $x,y \in B_E$ such that:
$$\displaystyle |z| \leq \frac{1+\max\{\|x\|,\|y\|\}}{2 \max\{\|x\|,\|y\|\}}.$$ 
\noindent Then $|1-p| \leq 2 |1-|z|^2 p|$ where $p$ denotes the scalar product $\langle x,y \rangle$.
\end{lema}

\dem  Suppose without loss of generality that $\|x\| \geq \|y\|$ and $x \neq 0$. Otherwise, the inequality is clearly true for any $z \in \C$. Notice that:
$$|1-p|=|1-|z|^2p+|z|^2p-p| \leq |1-|z|^2p|+||z|^2-1||p|,$$
\noindent so it is sufficient to prove that $|1-|z|^2p|+||z|^2-1||p| \leq 2|1-|z|^2p|$, which is equivalent to $|1-|z|^2p| \geq ||z|^2-1||p|$. We consider two cases: \vspace{0.1cm}

\vspace{-0.4cm}
i) if $|z|^2 \leq 1$: then $1-|z|^2 \geq 0$ so we need to prove $|1-|z|^2 p| \geq (1-|z|^2)|p|$ which is clearly satisfied since:
$$|1-|z|^2p| \geq 1-|z|^2 |p| \geq |p|-|z|^2 |p|=(1-|z|^2)|p|$$
\noindent where second inequality is satisfied since $|p| \leq \|x\| \|y\| < 1$. \smallskip

ii) On the other hand, suppose that $|z|^2 >1$: we need to prove that $|1-|z|^2p| \geq (|z|^2-1)|p|$. Since $|1-|z|^2p| \geq 1-|z|^2 |p|$, it is sufficient to prove that:
\ba \label{ineq1}
1-|z|^2 |p| \geq (|z|^2-1)|p|
\ea

\noindent Notice that $1-|z|^2|p| >0$ since $|z|^2 |p| \leq \|zx\| \|zy\| < 1$ because $zx,zy \in B_E$. Inequality (\ref{ineq1}) is equivalent to $2|z|^2 |p| <1+|p|$ which is true since:
$$2|z|^2|p| \leq 2\left( \frac{1+\|x\|}{2\|x\|}\right)^2 |p|$$
\noindent so we need to prove that: 
\begin{eqnarray*}
2\left( \frac{1+\|x\|}{2\|x\|}\right)^2 |p| \leq 1+|p| \mbox{ if and only if } 
\left(\frac{(1+\|x\|)^2}{2\|x\|^2}-1 \right)|p| \leq 1.
\end{eqnarray*}

\noindent But:
\begin{eqnarray*} 
\left(\frac{(1+\|x\|)^2}{2\|x\|^2}-1 \right)|p|=\frac{1+2\|x\|-\|x\|^2}{2\|x\|^2}|p| = \frac{1+(2-\|x\|)\|x\|}{2\|x\|^2}|p| \leq \\ \frac{1+1}{2\|x\|^2}\|x\|^2=1
\end{eqnarray*}
\noindent where last inequality is true because of the arithmetic mean-geometric mean inequality and since $|p| \leq \|x\| \|y\| \leq \|x\|^2$. \qed \bigskip



\smallskip

\subsubsection{\textbf{The case $E=\C$}}

If we deal with $E=\C$, it is easy to prove that (\ref{condi2}) is bounded:
\begin{prop} \label{clasic}
	Let $x,y \in \D$ and $z \in \C$ such that:
	$$|z| \leq \frac{1+\max\{|x|,|y|\}}{2 \max\{|x|,|y|\}}.$$
	\noindent Then $zx,zy \in \D$ and:
	$$\frac{\rho(zx,zy)}{|z|\rho(x,y)} \leq 2.$$
\end{prop}

\dem Suppose without loss of generality that $|x| \geq |y|$. Notice that $zx,zy \in \D$ since:
$$|zy| \leq |zx| \leq \frac{1+|x|}{2|x|} |x| = \frac{1+|x|}{2} < 1.$$
We have:
$$\rho(zx,zy)=\left| \frac{zx-zy}{1-\overline{zx}zy}\right|=|z|\left|  \frac{x-y}{1-|z|^2\overline{x}y}\right|=|z| \frac{|x-y|}{|1-|z|^2\overline{x}y|}$$ 
\noindent and:
$$|z| \rho(x,y)=|z| \left| \frac{x-y}{1-\overline{x}y}\right|=|z| \frac{|x-y|}{|1-\overline{x}y|},$$
\noindent so the inequality is equivalent to: 
$$|z| \frac{|x-y|}{|1-|z|^2\overline{x}y|} \leq 2 |z| \frac{|x-y|}{|1-\overline{x}y|}$$
\noindent or equivalently:
$$|1-\overline{x}y| \leq 2 |1-|z|^2\overline{x}y|.$$

\noindent Calling $p=\overline{x}y$, we have to prove that $|1-p| \leq 2 |1-|z|^2 p|$. Apply Lemma \ref{lema1} for $E=\C$ and we are done. \qed

\begin{rem}
Notice that the bound $2$ is the best possible. Indeed, take $x_n,y_n \in \D$ such that $x_n \rightarrow 1$ and $y_n \rightarrow -1$. It is clear that for $z_n \rightarrow 0$, the expression $|1-\overline{x_n}y_n|$ tends to $2$ when $n \rightarrow \infty$ and the expression $2 |1-|z_n|^2\overline{x_n}y_n|$ also tends to $2$, so the inequality above is sharp. \qed
\end{rem} 

\subsubsection{\textbf{The case when $E$ is any complex Hilbert space}}

In the following lemma, we will consider $x,y \in \B_E$ and $z \in \C$ such that $zx,zy \in B_E$ and we will denote by $r=\|x\|$, \ $s=\|y\|$ and $p = \langle x,y \rangle$. We will also denote $m=\Re p$ and $u=r^2+s^2$. Finally, we will also consider $A=\|x-y\|^2=r^2+s^2-2 \Re p=u-2m$ and $B=r^2s^2-|p|^2$. This notation will be also used in Theorem \ref{teo}. Notice that $A-B \geq 0$ since $A-B=r^2+s^2-2m-r^2s^2+|p|^2=|1-p|^2-(1-r^2)(1-s^2) \geq 0$ if and only if $|1-p|^2 \geq (1-r^2)(1-s^2)$ if and only if $1-\frac{(1-r^2)(1-s^2)}{|1-p|^2}=\rho(x,y)^2 \geq 0$ which is clearly true. 

\begin{lema} \label{lema2}
Let $E$ be a complex Hilbert space and $x,y \in B_E$. 
We have that:
$$\frac{\|x-y\|^2|1-p|^2}{\|x-y\|^2-(\|x\|^2 \|y\|^2-|p|^2)} \leq 4 \mbox{ or, equivalently: }\frac{A|1-p|^2}{A-B} \leq 4.$$
\end{lema}

\dem The inequality is equivalent to:
\begin{eqnarray*}
(1-2\Re p +|p|^2) \|x-y\|^2 \leq 4\|x-y\|^2 -4r^2s^2+4|p|^2 \mbox{ if and only if:} \\
(1-2\Re p +|p|^2) (r^2+s^2-2\Re p) \leq 4(r^2+s^2)-8 \Re p -4r^2s^2+4|p|^2.
\end{eqnarray*}
\noindent Bearing in mind that $m=\Re p$ and $u=r^2+s^2$, we need to prove:
\begin{eqnarray*}
(1-2m +|p|^2) (u-2m) \leq 4u-8m -4r^2s^2+4|p|^2 \leftrightarrow \\
u-2um+u |p|^2-2m+4m^2-2m|p|^2 \leq 4u-8m-4r^2s^2+4|p|^2 \leftrightarrow \\
4u-8m-4r^2s^2+4|p|^2-u+2um-u|p|^2+2m-4m^2+2m|p|^2 \geq 0.
\end{eqnarray*}

\noindent Notice that $|p| \geq |m|$ and since $4+2m-u \geq 0$, then $(4+2m-u) |p|^2 \geq (4+m-u) m^2$, so it is sufficient to prove:
$$4u-8m-4r^2s^2+(4+2m-u)m^2-u+2um+2m-4m^2 \geq 0.$$

It is also clear that $u \geq 2rs$, so $u^2 \geq 4r^2s^2$ and last inequality is equivalent to:
\begin{eqnarray*}
4u-8m-u^2+(4+2m-u)m^2-u+2um+2m-4m^2 \geq 0 \leftrightarrow \\
2m^3-u m^2+2(u-3)m+3u-u^2 \geq 0
\end{eqnarray*}
\noindent The expression at left can be easily factorized and equals to:
$$(u-2m)(3-u-m^2)$$
\noindent where both factors are clearly greater or equal to $0$ and we are done. \qed \bigskip

The following lemma will be used at the end of the proof ot Theorem \ref{teo}:
\begin{lema} \label{max_fun}
Let $f(a,b,c)=(3-b^2)(a-c)-(a^2-b^2)(2-c)$. Then $f(a,b,c) \geq 0$ for any $0 \leq c \leq b \leq a \leq 1$.
\end{lema}

\dem Notice that $f(a,b,c)=(3-b^2)a-2(a^2-b^2)-(3-a^2)c$, so $f$ is affine with respect to $c$. Hence it is enough to prove the inequality for $c=b$. The function becomes $f(a,b,b)=(3-b^2)(a-b)-(a^2-b^2)(2-b)=(a-b)((3-b^2)-(a+b)(2-b))$ and since $a-b \geq 0$, it is enough to prove that $(3-b^2)-(a+b)(2-b) \geq 0$.   The expression $g(a,b)=(3-b^2)-(a+b)(2-b)$ is affine with respect to $a$, so it is enough to prove it for $a=1$. Notice that $g(1,b)=(3-b^2)-(1+b)(2-b)=1-b$ which is clearly greater or equal to $0$, so we are done. \qed \bigskip

\begin{teo} \label{teo}
Let $E$ be a finite or infinite dimensional complex Hilbert space, $z \in \C$ and $x,y \in B_E$. If:
$$|z| \leq \frac{1+\max\{\|x\|,\|y\|\}}{2 \max\{\|x\|,\|y\|\}},$$
\noindent then $zx,zy \in B_E$ and:
$$\frac{\rho_E(zx,zy)}{|z|\rho_E(x,y)} \leq 2.$$
\end{teo}

\dem Suppose without loss of generality that $\|x\| \geq \|y\|$ and $z \neq 0$. We will denote $\rho=\rho_E(x,y)$ and $\rho_z=\rho_E(zx,zy)$. If $ \frac{1}{2} \leq |z| < 1$, then the result is clear since: 
$$\rho_E(zx,zy) \leq \rho_E(x,y) \leq 2 |z| \rho_E(x,y)$$
\noindent where first inequality is true because of the contractivity of the pseudohyperbolic distance for the function $g: B_E \to B_E$ given by $g(x)=z x$. 

So let us prove it for $|z| < 1/2$ or $|z| \geq 1$. Taking squares, the inequality is equivalent to prove:
$$\frac{\rho_z^2}{|z|^2\rho^2} \leq 4.$$

\noindent Bear in mind the expression (\ref{pseudoh}) for the pseudohyperbolic distance and call $t=|z|^2$. So we have:
\ban
\displaystyle  \frac{\rho_z^2}{|z|^2\rho^2}= \frac{\frac{|1-t p|^2-(1-t r^2)(1-ts^2)}{|1-tp|^2}}{\frac{t(|1-p|^2-(1-r^2)(1-s^2))}{|1-p|^2}} = \\ 
\frac{(|1-t p|^2-(1-t r^2)(1-ts^2))|1-p|^2}{t(|1-p|^2-(1-r^2)(1-s^2))|1-tp|^2}= \\
 \frac{(1+t^2|p|^2-2t\Re p-1-t^2 r^2s^2+t(r^2+s^2))|1-p|^2}{t(1+|p|^2-2\Re p-1-r^2s^2+r^2+s^2)|1-tp|^2}= \\
 \frac{(t^2|p|^2-2t\Re p-t^2 r^2s^2+t(r^2+s^2))|1-p|^2}{t(|p|^2-2\Re p-r^2s^2+r^2+s^2)|1-tp|^2} 
\ean

\noindent and dividing by $t$:
\ban
\frac{\rho_z^2}{|z|^2\rho^2} \leq  \frac{(t|p|^2-2\Re p-t r^2s^2+r^2+s^2)|1-p|^2}{(|p|^2-2\Re p-r^2s^2+r^2+s^2)|1-tp|^2}= \\  \frac{(\|x-y\|^2-t(r^2s^2-|p|^2))|1-p|^2}{(\|x-y\|^2-(r^2s^2-|p|^2))|1-tp|^2}.
\ean

\noindent Using the notation above, last inequality is equivalent to:
\ba \label{main_ineq}
\frac{(A-tB) |1-p|^2}{(A-B)|1-tp|^2} \leq 4
\ea 
\noindent so we need to prove (\ref{main_ineq}) for $t \geq 1$ and $t \leq 1/4$.

\noindent If $t \geq 1$, then the result is clear since:
\begin{eqnarray*}
\frac{(A-tB) |1-p|^2}{(A-B)|1-tp|^2} \leq  \left(\frac{A-B}{A-B}\right) \frac{|1-p|^2}{|1-tp|^2} = \frac{|1-p|^2}{|1-tp|^2} \leq 4
\end{eqnarray*}
where last inequality is true by Lemma \ref{lema1}. So we need to prove inequality (\ref{main_ineq}) for $0 \leq t \leq 1/4$. This inequality is clearly equivalent to:
\begin{eqnarray} \label{ineq21}
4 (A-B)|1-tp|^2 \geq (A-tB)|1-p|^2
\end{eqnarray}
\noindent which is satisfied if and only if:
\begin{eqnarray*}
4 (A-B)(1-2mt+t^2|p|^2) \geq (A-tB)|1-p|^2.
\end{eqnarray*}
\noindent This is equivalent to:
\begin{eqnarray*}
4 (A-B)|p|^2t^2 +(B|1-p|^2-8m(A-B))t+4(A-B)-A|1-p|^2 \geq 0. 
\end{eqnarray*}

Since $B,t \geq 0$ and $4(A-B)-A|1-p|^2 \geq 0$ by Lemma \ref{lema2}, the inequality is clearly true if $m <0$. So we can suppose without loss of generality that $m \geq 0$. We will prove inequality (\ref{ineq21}) for $m \geq 0$. The inequality is equivalent to:
\begin{eqnarray*}
4 (A-B)|1-tp|^2 - (A-tB)|1-p|^2 \geq 0
\end{eqnarray*}
so we will prove last inequality. Notice that:
\begin{eqnarray*}
4 (A-B)|1-tp|^2 - (A-tB)|1-p|^2 =\\
4(A-B)(1-2mt+t^2 |p|^2)-(A-B)(1-2m+|p|^2)-B(1-t)|1-p|^2=\\
4(A-B)(1-2m+|p|^2+2m(1-t)-(1-t^2)|p|^2)\\-(A-B)(1-2m+|p|^2)-B(1-t)|1-p|^2=\\
3(A-B)(1-2m+|p|^2)+8m(1-t)(A-B)\\-4(1-t^2)|p|^2(A-B)-B(1-t)|1-p|^2
\end{eqnarray*}
\noindent Since $0 \leq t \leq 1/4$, we have that $3/4 \leq 1-t \leq 1$. Hence:
\begin{eqnarray*}
4 (A-B)|1-tp|^2 - (A-tB)|1-p|^2 \geq \\
3(A-B)(1-2m+|p|^2)+8m(1-t)(A-B)\\-4(1-t^2)|p|^2(A-B)-B(1-t)|1-p|^2 \geq \\
3(A-B)-6m(A-B)+3|p|^2(A-B)+8m \cdot \frac{3}{4}(A-B)\\-4|p|^2(A-B)-B|1-p|^2 =\\
3(A-B)-6m(A-B)+6m(A-B)-|p|^2(A-B)-B|1-p|^2 =\\
(3-|p|^2)(A-B)-B|1-p|^2.
\end{eqnarray*}

Notice that:
\begin{eqnarray*}
(3-|p|^2)(A-B)-B|1-p|^2 = \\ (3-|p|^2)A-B(3-|p|^2+1-2m+|p|^2)=\\
(3-|p|^2)A-B(4-2m) = \\ (3-|p|^2)(r^2+s^2-2m)-(4-2m)(r^2s^2-|p|^2) \geq  \\ (3-|p|^2)(2rs-2m)-(4-2m)(r^2s^2-|p|^2)= \\ 2(3-|p|^2)(rs-m)-2(2-m)(r^2s^2-|p|^2)= \\ 2(3-b^2)(a-c)-2(a^2-b^2)(2-c).
\end{eqnarray*}

To finish the proof, notice that if we call $a=rs$, $b=|p|$ and $c=m$, then $0 \leq c \leq b \leq a \leq 1$ so we use Lemma \ref{max_fun} and:
\begin{eqnarray*}
(3-|p|^2)(A-B)-B|1-p|^2 \geq \\ 2(3-b^2)(a-c)-2(a^2-b^2)(2-c)=2 f(a,b,c) \geq 0
\end{eqnarray*} 
\noindent and we are done. \qed \bigskip

\subsubsection{Results for $X=C_0(S)$} \label{resultadoscox}

Let $S$ be a locally compact topological space and consider $X=C_0(S)$ given by the space of continuous functions $f: S \to \C$  such that for any $\eps >0$, there exists a closed compact subset $K \subset S$ such that $|f(x)| < \eps$ for any $x \in S \setminus K$. Endowed with the sup-norm, $C_0(S)$ becomes a Banach space and the pseudohyperblic distance for $x,y \in C_0(S)$ is well-known (see \cite{AGL}) and it is given by:
\ba \label{pseudo_cox}
\rho_X(x,y)=\sup_{t \in S} \rho(x(t),y(t)).
\ea 

We prove that expression (\ref{condi2}) is also bounded by $2$ when we deal with the space $X=C_0(S)$:

\begin{prop}
Let $X=C_0(S)$ and $x,y \in X$. If $z \in \C$ satisfies:
$$|z| \leq \frac{1+\max\{\|x\|,\|y\|\}}{2 \max \{\|x\|,\|y\|\}},$$
\noindent then:
$$\frac{\rho_X(zx,zy)}{|z| \rho_X(x,y)} \leq 2.$$
\end{prop}

\dem Suppose without loss of generality that $\|x\| \geq \|y\|$. For any $t \in S$, we have that $x(t),y(t) \in \D$ since $\|x\|=\sup_{t \in S}|x(t)| <1$ and $\|y\|=\sup_{t \in S} |y(t)|<1$. The result is clear since:
\begin{eqnarray*}
\rho_X(zx,zy)=\sup_{t \in S} \rho(z x(t),zy(t)) \leq \sup_{t \in S} 2|z| \rho(x(t),y(t))=\\ 2|z| \sup_{t \in S}  \rho(x(t),y(t))=2 |z| \rho_X(x,y)
\end{eqnarray*}
\noindent where first inequality is clear because of Proposition \ref{clasic} and because for any $t \in X$ we have that:
$$|z| \leq \frac{1+\|x\|}{2\|x\|}=\frac{\frac{1}{\|x\|}+1}{2} \leq \inf_{t \in X} \left\{\frac{1+|x(t)|}{2|x(t)|}\right\} \leq \frac{1+|x(t)|}{2|x(t)|}$$
\noindent and we are done. \qed \bigskip

\section{Applications} \label{section4}

Now we give some applications related to Theorem \ref{teo}. First, we will show that the function $x \mapsto (1-\|x\|^2) |\R f(x)|$ for $x \in B_E$ is Lipschitz with respect to the pseudohyperbolic distance. Hence, we derive some results about interpolating sequences for $\B(B_E)$ in Subsection \ref{subsec42}. Indeed, we provide a new proof that these sequences are separated for the pseudohyperbolic distance. We also prove that these sequences can be slightly perturbed and they remain interpolating. Finally, in Subsection \ref{subsec43} we will take an in-depth look at the automorphisms of $B_E$. This will permit us to give necessary and sufficient conditions for a composition operator on $\B(B_E)$ to be bounded below. 

\subsection{The Lipschitz continuity of $(1-\|x\|^2) |\R f(x)|$} \label{subsec31}

We will denote by $\Pi$ the unit circle of the complex plane $\C$, that is, the set of complex numbers $u$ such that $|u|=1$.

\begin{lema} \label{lema_eps}
	Let $f \in \B(B_E)$. Fix $\eps >0$ and $x,y \in B_E$. If $(1+\eps u)x$ and $(1+\eps u)y$ belongs to $B_E$ for any $u \in \Pi$, then there exists $u_0 \in \Pi$ such that:
	\ba
	|\R f(x)-\R f(y)| \leq \frac{1}{\eps} \|f\|_\I \beta((1+\eps u_0) x,(1+\eps u_0)y).
	\ea
\end{lema}

\dem Fix $x,y \in B_E$ and $\eps >0$. Notice that the function $f(x+\eps u x)-f(y+\eps u y)$ defined for $u \in \Pi$ is continuous. Since $\Pi$ is a compact set, there exists $u_0 \in \Pi$ such that:
$$f(x+\eps u_0 x)-f(y+\eps u_0 y)= \max \{f(x+\eps u x)-f(y+\eps u y) : u \in \Pi \}.$$

\noindent Consider $g(u)=f(x+\eps u x)$ for $u$ defined on an open disk of the complex plane $\C$ which contains $\Pi$. It is clear that:
$$g'(u)=\nabla f(x+\eps u x)(\eps x),$$

\noindent so $g'(0)=\eps \R f(x)$. Similarly, if $h(u)=f(y+\eps u y)$, then $h'(0)=\eps R f(y)$. By the Cauchy's integral formula we have:
\ban
|\R f(x)-\R f(y)|=|\langle x,\overline{\nabla f(x)} \rangle -\langle y,\overline{\nabla f(y)} \rangle|=\\ \left| \frac{1}{\eps} \frac{1}{2 \pi i} \int_{|u|=1} f(x+\eps u x)-f(y+\eps u y) \frac{du}{u^2}\right| \leq \\ \frac{2\pi}{2 \pi \eps} |f(x+\eps u_0 x)-f(y+\eps u_0 y) | \leq \frac{1}{\eps} \|f\|_\I \beta((1+\eps u_0) x,(1+\eps u_0)y)
\ean
\noindent where last inequality is true by Proposition \ref{prop_lip}. \qed \bigskip

The proof of the following lemma is an easy calculation. It will be used in Lemma \ref{lema43}.
\begin{lema} \label{desig_log}
For any $0 \leq t < 1$ we have:
$$\frac{1}{2} \log \left(\frac{1+t}{1-t}\right) \leq \frac{t}{1-t}.$$
\end{lema}

\begin{lema} \label{lema43}
	Let $f \in \B(B_E)$ and $x,y \in B_E$ such that $\|x\| \geq \|y\|$. Then:
	\ban
	(1-\|x\|^2) |\R f(x)-\R f(y)| \leq 12 \|f\|_\I  \rho_E(x,y).
	\ean
\end{lema}

\dem Take:
$$\eps=\frac{1-\|x\|}{2\|x\|} >0.$$
\noindent Notice that for any $u \in \Pi$ we have that $(1+\eps u)x$ and $(1+\eps u) y$ belongs to $B_E$ since:
$$(1+\eps) \|x\| \leq \left( 1+\frac{1-\|x\|}{2\|x\|} \right) \|x\|=\frac{1+\|x\|}{2 \|x\|} \|x\| = \frac{1+\|x\|}{2} < 1$$
\noindent so clearly $\|(1+\eps u) x\| \leq (1+\eps) \|x\| < 1$
\noindent and since $\|y\| \leq \|x\|$:
$$\|(1+\eps u) y\| \leq (1+\eps) \|y\| \leq (1+\eps) \|x\| < 1.$$

\noindent By Lemma \ref{lema_eps}, there exists $u_0 \in \Pi$ such that:
	\ban
|\R f(x)-\R f(y)| \leq \frac{2\|x\|}{1-\|x\|} \|f\|_\I \beta_E((1+\eps u_0) x,(1+\eps u_0)y).
\ean

 Take $z_0=1+\eps u_0$ which satifies:
$$|z_0| \leq 1+\eps =1+\frac{1-\|x\|}{2\|x\|}= \frac{1+\|x\|}{2\|x\|}.$$


\noindent By Theorem \ref{teo} we have that $z_0 x,z_0 y \in B_E$ and:
\ba \label{desig11}
\rho_E(z_0 x,z_0y) \leq 2 |z_0| \rho_E(x,y).
\ea

\noindent Denote $B=(1-\|x\|^2) |\R f(x)-\R f(y)|$. We obtain:
$$B \leq (1-\|x\|^2) \frac{2\|x\|}{1-\|x\|} \|f\|_\I \beta_E(z_0x,z_0y).$$


By Lemma \ref{desig_log} we also have:
$$\beta_E (z_0 x,z_0y) \leq \frac{\rho_E(z_0x,z_0y)}{1-\rho_E(z_0 x,z_0y)}$$
\noindent  We obtain:
\ban
B \leq (1+\|x\|)(1-\|x\|) \frac{2\|x\|}{1-\|x\|} \|f\|_\I \frac{\rho_E(z_0x,z_0y)}{1-\rho_E(z_0 x,z_0y)} \leq \\
4 \|x\| \|f\|_\I \frac{\rho_E(z_0x,z_0y)}{1-\rho_E(z_0 x,z_0y)} = \frac{4 \|x\| \|f\|_\I}{\frac{1}{\rho_E(z_0 x,z_0y)}-1}
\ean

\noindent so:
$$\left( \frac{1}{\rho_E(z_0 x,z_0y)}-1\right) B \leq 4 \|x\| \|f\|_\I$$
\noindent which is equivalent to:
\ba \label{desig12}
\frac{B}{\rho_E(z_0 x,z_0y)} -B \leq 4 \|x\| \|f\|_\I.
\ea

\noindent Bear in mind (see (\ref{seminormas})) that $\|f\|_\B \leq \|f\|_\I$. We have:
\ban
B \leq (1-\|x\|^2) |\R f(x)|+(1-\|x\|^2) |\R f(y)| \leq \\ (1-\|x\|^2) \| \nabla f(x) \| \|x\| +(1-\|y\|^2) \| \nabla f(y)\| \|y\| \leq \\ 2\|f\|_\B \|x\| \leq 2  \|x\|\|f\|_\I
\ean
\noindent so from inequality (\ref{desig12}) we have:
$$\frac{B}{\rho_E(z_0 x,z_0y)} \leq 4 \|x\| \|f\|_\I + B \leq 4 \|x\|\|f\|_\I + 2 \|x\|\|f\|_\I = 6\|x\| \|f\|_\I$$
\noindent and we conclude $B \leq 6 \|x\| \|f\|_\I \rho_E(z_0 x,z_0y)$. Finally, we apply inequality (\ref{desig11}) and since $|z_0| \leq \frac{1+\|x\|}{2\|x\|}$ we obtain:
\ban
B \leq 6 \|x\| \|f\|_\I 2 |z_0| \rho_E(x,y) \leq 12 \|x\| \|f\|_\I \frac{1+\|x\|}{2\|x\|} \rho_E(x,y) = \\ 12 \|f\|_\I \frac{1+\|x\|}{2} \rho_E(x,y) \leq 12 \|f\|_\I \rho_E(x,y)
\ean
\noindent and we are done. \qed \bigskip

\begin{teo} \label{teo1}
	Let $f \in \B(B_E)$ and $x,y \in B_E$. Then:
	$$|(1-\|x\|^2) \R f(x)-(1-\|y\|^2) \R f(y)| \leq 14 \|f\|_\I \rho_E(x,y).$$
\end{teo}

\dem Call $F=|(1-\|x\|^2) \R f(x)-(1-\|y\|^2) \R f(y)|$ and suppose without loss of generality that $\|x\| \geq \|y\|$. We have that:
\ba \label{desig13} \
F=|(1-\|x\|^2) (\R f(x)-\R f(y))-(\|x\|^2-\|y\|^2) \R f(y)| \leq \nonumber \\
(1-\|x\|^2) |\R f(x)-\R f(y)| +(\|x\|^2-\|y\|^2) |\R f(y)|.
\ea

\noindent Since $\|x\|^2-\|y\|^2=(\|x\|+\|y\|)(\|x\|-\|y\|) \leq 2 (\|x\|-\|y\|)$ and bearing in mind that $\rho_E(\|x\|,\|y\|) \leq \rho_E(x,y)$ we obtain:
\ban 
(\|x\|^2-\|y\|^2) |\R f(y)| \leq 2\frac{\|x\|-\|y\|}{1-\|x\| \|y\|} (1-\|x\| \|y\|) |\R f(y)| \leq \\ 2 \rho_E(\|x\|,\|y\|)(1-\|y\|^2)|\R f(y)| \leq 2 \|y\|\|f\|_\B \rho_E(x,y) \leq  2 \|f\|_\I \rho_E(x,y).
\ean

By Lemma \ref{lema43} we know that:
$$(1-\|x\|^2) |\R f(x)-\R f(y)| \leq 12 \|f\|_\I \rho_E(x,y)$$
\noindent so from (\ref{desig13}) we conclude:
$$F \leq 12 \|f\|_\I \rho_E(x,y)+2\|f\|_\I \rho_E(x,y)=14 \|f\|_\I \rho_E(x,y)$$
\noindent and we are done. \qed \bigskip

The following corollary extends to the finite and infinite dimensional setting results given by Attele in \cite{At} and improved by Xiong in \cite{Xi} for the classical Bloch space $\B$.

\begin{cor} \label{cor45}
	Let $E$ be a complex Hilbert space. The function $x \mapsto (1-\|x\|^2) |\R f(x)|$ for $x \in B_E$ is Lipschitz with respect to the pseudohyperbolic distance and the following inequality holds:
	$$|(1-\|x\|^2) |\R f(x)|-(1-\|y\|^2) |\R f(y)| | \leq 14 \|f\|_\I \rho_E(x,y).$$
\end{cor}

\dem Applying Theorem \ref{teo1} it is clear that:
\ban
|(1-\|x\|^2) |\R f(x)|-(1-\|y\|^2) |\R f(y)| | \leq \\
|(1-\|x\|^2) \R f(x)-(1-\|y\|^2) \R f(y)| \leq 14 \|f\|_\I \rho_E(x,y).
\ean
and we are done. \qed \bigskip

\subsection{Results on interpolating sequences for the Bloch space} \label{subsec42}

Let $E$ be a complex finite or infinite dimensional complex Hilbert space. Recall that $(x_n) \subset B_E \setminus \{0\}$ is said to be interpolating for the Bloch space $\B(B_E)$ if for any bounded sequence $(a_n)$ of complex numbers, there exists $f \in \B(B_E)$ such that $(1-\|x_n\|^2) \R f(x_n)=a_n$. Attele studied in \cite{At} this kind of interpolation for the classical Bloch space $\B$ and the finite and infinite dimensional setting was studied in \cite{BGLM2}. We provide a new approach to prove that a necessary condition for a sequence $(x_n) \subset B_E$ to be interpolating for $\B(B_E)$ is to be separated for the pseudohyperbolic distance:

\begin{prop}
Let $E$ be a complex Hilbert space.	If $(x_n) \subset B_E \setminus \{0\}$ is interpolating for $\B(B_E)$, then there exists $C>0$ such that $\rho(x_k,x_j) \geq C$ for any $k \neq j$, $k,j \in \N$.
\end{prop}

\dem Since $(x_n) \subset B_E \setminus \{0\}$ is interpolating, there exists a sequence $(f_n) \subset \B(B_E)$ such that:
\ban (1-\|x_n\|^2) \R f_n(x_n) =1 \ \mbox{ and } \ (1-\|x_k\|^2) \R f_n(x_k) =0 \ \mbox{ if } \ k \neq n.
\ean
\noindent The operator $T: \B(B_E) \to \ell_{\infty}$ given by $T(f)=((1-\|x_n\|^2) \R f(x_n))$ is surjective, so by the Open Mapping Theorem, there exists $M >0$ such that $\|f\|_R \leq M \sup_{j \in \N} (1-\|x_j\|^2) |\R f(x_j)|$ so $\|f_n\|_R \leq M$ for any $n \in \N$. Applying Theorem \ref{teo1}, we have:
\ban
|(1-\|x_n\|^2) \R f_n(x_n)-(1-\|x_k\|^2) \R f_n(x_k)| \leq 14 \|f_n\|_\I \rho_E(x_k,x_j) \leq \\ 14 A_0 \|f_n\|_\R \rho_E(x_k,x_j) \leq 14 A_0  M \rho_E(x_k,x_j). 
\ean
\noindent Hence, $1-0 \leq 14 A_0 M \rho_E(x_k,x_j)$ and we conclude that:
$$\rho_E(x_k,x_j) \geq \frac{1}{14 A_0 M}$$
\noindent so we are done. \qed \bigskip

Attele (see \cite{At}) also proved that any interpolating sequence $(z_n) \subset \D$ for $\B$ can be slightly perturbed and the sequence remains interpolating. By means of Theorem \ref{teo1}, we adapt his proof and generalize the result to the case when we deal with any complex Hilbert space $E$.
\begin{teo}
	If $(x_n) \subset B_E \setminus \{ 0 \}$ is an interpolating sequence for $\B(B_E)$, then there exists $\delta >0$ such that if $(y_n) \subset B_E $ satisfies that $\sup_{n \in \N} \rho_E(x_n,y_n)<\delta$ then $(y_n)$ is also an interpolating sequence for $\B(B_E)$.
\end{teo}

\dem Since $(x_n)$ is interpolating, the operator $T: \B(B_E) \to \ell_\infty$ given by $T(f)=((1-\|x_n\|^2) \R f(x_n))$ is surjective. Hence, its adjoint $T^*: \ell_{\infty}^* \to (\B(B_E))^*$ is injective and it has closed range. In particular, $T^*$ is left-invertible. The set of left-invertible elements is open in the Banach algebra of linear operators from $\ell_{\infty}^*$ to $(\B(B_E))^*$. So there exists $\delta$ such that if $\|T^*-R\| < 14 A_0 \delta$, then $R$ is left-invertible. If we consider $S(f)=((1-\|y_n\|^2) \R f(y_n))$, then by Theorem \ref{teo1}:
\ban
\|(T-S)(f)\|_{\infty} = \sup_{n \in \N} |(1-\|x_n\|^2) \R f(x_n) -(1-\|y_n\|^2) \R f(y_n)|\leq \\ 14 \rho_E(x_n,y_n) \|f\|_\I \leq 14 A_0 \|f\|_\R \rho_E(x_n,y_n) < 14 A_0 \delta \|f\|_\R
\ean

\noindent so $\|T-S\| < 14 A_0 \delta$ and hence $\|T^*-S^*\|=\|T-S\| < 14 A_0 \delta$. We conclude that $S^*$ is left-invertible and hence $S$ is surjective, as we wanted. \qed \bigskip

\subsection{Bounded below composition operators} \label{subsec43}

Recall that a linear operator between Banach spaces $T: X \to Y$ is bounded below if there exists $k >0$ such that $\|x\| \leq k \|T(x)\|$. It is well-known that the linear bounded operator $T$ is bounded below if and only if $T$ is injective and has closed range. \medskip

Let $\phi: \D \to \D$ be an analytic map.  The composition operator $C_{\phi}: \B \to \B$ is given by $C_{\phi}(f)=f \circ \phi$. It is well-known that $C_{\phi}$ is bounded for any $\phi$. Denote by:
\ba \label{tau1}
\tau_{\phi}(z)=\frac{1-|z|^2}{1-|\phi(z)|^2} \phi'(z).
\ea
\noindent In \cite{GYZ}, the authors investigated conditions under which $\phi$ induces a composition operator with closed range on the Bloch space $\B$. In particular, they proved the following necessary condition:
\begin{prop}
	If $C_{\phi}$ is bounded below, then there exist $\eps, r >0$ with $r < 1$ such that for any $z \in \D$ we have $\rho(\phi(w),z) \leq r$ for all $w \in \D$ satisfying $|\tau_{\phi}(w)| > \eps$.
\end{prop}

\noindent In order to provide sufficient conditions, the authors studied the function given by $z \mapsto (1-|z|^2)|f'(z)|$ for $z \in \D$ and they proved this function is Lipschitz with respect to the pseudohyperbolic distance if $f \in \B$. Indeed, they proved:
\ba \label{lips}
\ \ |(1-|z|^2)|f'(z)|-(1-|w|^2)|f'(w)|| \leq 3.31 \|f\|_\B \rho(z,w)
\ea
\noindent for any $z,w \in \D$ and $f \in \B$. This result improves a previous result given by Attele in \cite{At} since the constant given by him was $9$ instead of $3.31$. Xiong improved the constant in \cite{Xi}, providing $3 \sqrt{3}/2 \approx 2.6$. Inequality (\ref{lips}) permitted to give the following sufficient condition for $C_{\phi}$ to be bounded below (see \cite{GYZ}):
\begin{teo} \label{teo5}
Let $\phi: \D \to \D$ an analytic map. Suppose that there exist $0 < r < \frac{1}{4}$ and $\eps >0$ such that for any $w \in \D$ there exists $z_w \in \D$ satisfying $\rho(\phi(z_w),w) < r$ and $|\tau_{\phi}(z_w)| > \eps$. Then, $C_\phi: \B \to \B$ is bounded below. 
\end{teo}

Some authors (see \cite{Ch} and \cite{DJO}) extended these results by considering analytic maps $\phi: B_n \to B_n$, where $B_n$ denotes the open unit ball in the finite dimensional Hilbert space $(\C^n,\| \cdot \|_2)$. Nevertheless, the authors substituted $\tau_\phi(z)$ given in (\ref{tau1}) by the expression:
\ba \label{tau2}
\left( \frac{1-\|z\|^2}{1-\|\phi(z)\|^2}\right)^{(n+1)/2} |det(J_{\phi}(z))|
\ea
\noindent where $J_{\phi}(z)$ denotes the Jacobian $n \times n$ matrix of $\phi$. Indeed, $\tau_\phi(z)=1$ if $\phi$ is an automorphism of $B_n$. Furthermore, the authors also based their proofs on the definition of Bloch function on $B_n$ introduced by Timoney (see \cite{T80}) using the Bergman metric. \medskip

We want to extend the necessary the classical results on $\B$ to the finite and infinite dimensional setting, so we will provide necessary and sufficient conditions avoiding expression (\ref{tau2}) or the use of the Bergman metric. So, consider a complex finite or infinite dimensional Hilbert space $E$ and let $\psi: B_E \to B_E$ an analytic map. In order to extend the classical results, we introduce $\tau_\psi(x)$ and $\widetilde{\tau_\psi}(x)$ for $x \in B_E$ given by:
\ba \label{tau4}
{\tau_\psi}(x) =\frac{1-\|x\|^2}{1-\|\psi(x)\|^2} \|\psi'(x)\|.
\ea
\noindent and:
\ba \label{tau3}
\widetilde{\tau_\psi}(x)=\frac{\sqrt{1-\|x\|^2}}{1-\|\psi(x)\|^2} \|\psi'(x)\|
\ea
\noindent for $x \in B_E$. Notice that $\widetilde{\tau_\psi}(x) \geq {\tau_\psi}(x)$ for any $x \in B_E$.

\medskip

The boundness and compactness of the composition operator $C_{\psi}: \B(B_E) \to \B(B_E)$ given by $C_{\psi}(f)=f\circ \phi$ was studied in \cite{BGLM}. In particular, the authors proved that $C_{\psi}$ is bounded for any analytic map $\psi: B_E \to B_E$. In addition, the authors proved that $\|f \circ \psi\|_\I \leq \|f\|_\I$. \medskip

The following result will be used in Lemma \ref{lema_simp}:

\begin{lema} \label{lema_calc}
	Let $E$ be a complex Hilbert space and $f \in \B(B_E)$. Then for any $x \in B_E$:
	$$|f(x)-f(0)| \leq \|x\| \frac{\|f\|_\B}{1-\|x\|^2}.$$
\end{lema}
\dem We have that:
\ban
|f(x)-f(0)| =\left| \left( \int_{0}^{1} f'(xt) dt \right)(x) \right| \leq \|x\| \left\| \int_{0}^{1} \frac{f'(xt) (1-\|tx\|^2)}{1-\|tx\|^2} dt \right\| \leq \\ \|x\| \|f\|_\B \int_{0}^{1} \left| \frac{1}{1-\|tx\|^2} \right| dt  \leq  \|x\| \|f\|_\B \int_{0}^{1} \frac{1}{1-\|x\|^2} dt = \|x\| \frac{\|f\|_\B}{1-\|x\|^2}
\ean
\noindent and we are done. \qed \bigskip

Since the semi-norms $\| \cdot \|_\R$, $\| \cdot \|_\I$ and $\| \cdot \|_\B$ and their corresponding norms $\| \cdot \|_{\R-Bloch}$, $\| \cdot \|_{\I-Bloch}$ and $\| \cdot \|_{\B-Bloch}$ are equivalent, we can consider any of them in order to study if $C_\psi: \B(B_E) \to \B(B_E)$ is bounded below.

\begin{lema} \label{lema_simp}
Let $E$ be a complex Hilbert space and let $\psi: B_E \to B_E$ be an analytic map. The composition operator $C_{\psi} : \B(B_E) \to \B(B_E)$ is bounded below if and only if there exists $k >0$ such that:
$$\|C_\psi(f)\|_\I \geq k \|f\|_\I.$$
\end{lema}

\dem Suppose that $C_{\psi}$ is bounded below and let $f \in \B(B_E)$. There exists $k >0$ such that $\|C_{\psi}(f)\|_{\I-Bloch} \geq k \|f\|_{\I-Bloch}$. Consider $g(x)=f(x)-f(\psi(0))$. It is clear that $g(\psi(0))=0$ so:
\ban
\|C_{\psi}(f)\|_{\I}= \|f \circ \psi \|_\I = \| g \circ \psi\|_\I = \| g \circ \psi\|_{\I-Bloch} \geq \\ k \|g\|_{\I-Bloch} \geq k \|g\|_\I =k \|f\|_\I.
\ean

\noindent Now supose that $\|C_\psi(f)\|_\I \geq k \|f\|_\I$ for $0 < k \leq 1$. We will prove that there exists $k' >0$ such that  $\|C_\psi(f)\|_{\I-Bloch} \geq k' \|f\|_{\I-Bloch}.$ By Lemma \ref{lema_calc} we have:
$$|f(\psi(0))-f(0)| \leq \|\psi(0)\| \frac{\|f\|_\B}{1-\|\psi(0)\|^2}$$
\noindent so:
\ban
|f(\psi(0))| \geq |f(0)|-\|\psi(0)\| \frac{\|f\|_\B}{1-\|\psi(0)\|^2} \geq |f(0)|-\frac{\|f\|_\I}{1-\|\psi(0)\|^2}.
\ean

\noindent and we get:
$$|f(\psi(0))|+\frac{1}{(1-\|\psi(0)\|^2)} \|f\|_\I \geq |f(0)|.$$

\noindent Hence:
\ban
k(1-\| \psi(0)\|^2)|f(\psi(0))|+\|C_{\psi}(f)\|_\I \geq \\ k(1-\| \psi(0)\|^2) |f(\psi(0))|+k \|f\|_\I \geq k (1-\| \psi(0)\|^2) |f(0)|
\ean
\noindent so we conclude:
\ban
2(|f(\psi(0))|+\|C_\psi (f)\|_\I) = 2|f(\psi(0))| + \|C_\psi(f)\|_\I +\|C_\psi(f)\|_\I \geq \\ k(1-\|\psi(0)\|^2)|f(\psi(0))|+ \|C_\psi(f)\|_\I +\|C_\psi(f)\|_\I \geq \\ k(1-\|\psi(0)\|^2)|f(0)|+\|C_\psi(f)\|_\I \geq k(1-\|\psi(0)\|^2)(|f(0)|+\|C_\psi(f)\|_\I)
\ean 
\noindent and we conclude:
\ban
\|C_\psi(f))\|_{\I-Bloch} \geq \frac{k(1-\|\psi(0)\|^2)}{2} \|f\|_{\I-Bloch}
\ean
so we take $k'=k(1-\|\psi(0)\|^2)/2$ and we conclude that $C_\psi$ is bounded below as we wanted. \qed \bigskip


\subsubsection{\textbf{Study of the automorphisms $\phi_x$}} \label{subsec331}

In order to study necessary and sufficient conditions for bounded below composition operators, we need to provide several calculations related to the automorphisms $\phi_x$ of $B_E$ introduced in (\ref{automorph}). It is well-known that if $E$ is a finite dimensional Hilbert space, then $\phi_x$ is an involution (see Theorem 2.2.2 in \cite{R4}). Nevertheless, the proof of this result makes use of the Cartan's uniqueness theorem. So our first result provides a new proof of this assertion without the use of the Cartan's theorem, that is, we prove that for any $x \in B_E$, the automorphism $\phi_x: B_E \to B_E$ is an involution for any finite or infinite dimensional complex Hilbert space $E$:
\begin{lema} \label{involution}
Let $E$ be a complex Hilbert space and $x \in B_E$. Then $\phi_x$ is an involution, that is, $\phi_x \circ \phi_x=Id_E$.
\end{lema}

\dem By (\ref{automorph}), we know that:
\ban
\phi_x(\phi_x(y))=(s_x Q_x+P_x)(m_x(\phi_x(y))= (s_x Q_x+P_x) \left( \frac{x-\phi_x(y)}{1-\langle \phi_x(y),x \rangle}\right) 
\ean
\noindent and bearing in mind (see Lemma 3.6 in \cite{Mir1}) that:
$$1-\langle \phi_x(y),x\rangle=1-\langle \phi_x(y),\phi_x(0)\rangle=\frac{1-\|x\|^2}{1-\langle y,x \rangle}$$
\noindent then:
\ban
\phi_x(\phi_x(y))=\frac{1-\langle y,x \rangle}{1-\|x\|^2} (s_x Q_x+P_x) \left(x-\phi_x(y)\right)=\\ \frac{1-\langle y,x \rangle}{1-\|x\|^2} \left((s_x Q_x+P_x)(x)-(s_x Q_x +P_x)((s_x Q_x+P_x)(m_x(y)))\right).
\ean
\noindent Since $P_x$ is an orthogonal projection and $Q_x$ is its orthogonal complement, we have that $P_x \circ Q_x=Q_x \circ P_x=0$ and $P_x+Q_x=Id_E$ and $P_x^2=P_x$, $Q_x^2=Q_x$ so:
\ban
\phi_x(\phi_x(y))=\frac{1-\langle y,x \rangle}{1-\|x\|^2} \left( x-(s_x^2 Q_x +P_x)\left( \frac{x-y}{1-\langle y,x \rangle} \right)\right)=\\ \frac{1-\langle y,x \rangle}{(1-\|x\|^2)(1-\langle y,x \rangle)} \left( (1-\langle y,x \rangle) x-(s_x^2 Q_x +P_x)\left( x-y \right)\right)= \\ \frac{1}{(1-\|x\|^2)} \left( (x-\|x\|^2 P_x(y)-x+(1-\|x\|^2) Q_x(y)+P_x(y)\right)=  \\ \frac{1}{(1-\|x\|^2)}(1-\|x\|^2)(P_x(y)+Q_x(y))=y
\ean
\noindent and we are done. \qed \bigskip

\begin{lema} \label{derivadas}
If $x \in B_E$ then $\phi_x'(0)$ is an invertible operator and $\phi_x'(0)^{-1}=\phi_x'(x)$.
\end{lema}

\dem By Lemma \ref{involution}, it is clear that $(\phi_x \circ \phi_x)'(0)=Id_E'(0)=Id_E$ so:
\ban
\phi_x'(\phi_x(0)) \circ \phi_x'(0)=\phi_x'(x) \circ \phi_x'(0)=Id_E
\ean
\noindent and we conclude the result. \qed \bigskip

As we mentioned in (\ref{invar}), $\|f\|_\I=\sup_{x \in B_E} \| \widetilde{\nabla} f(x) \|$. Notice that for any $x \in B_E$ we have:
\begin{align} \label{310}
\| \widetilde{\nabla} f(x) \|= \sup_{u \in \overline{B_E}} \| f'(\phi_x(0)) \circ \phi_x'(0)(u)\|=\sup_{w \in E \setminus \{0\}} \frac{| f'(x)(w)|}{\|\phi_x'(0)^{-1}(w)\|}
\end{align}
\noindent and for any $w \in E$ the expression $\|\phi_x'(0)^{-1}(w)\|^2$ is given (see \cite{BGM}) by:
\ba \label{coci}
\|\phi_x'(0)^{-1}(w)\|^2=\frac{(1-\|x\|^2) \|w\|^2+| \langle w,x \rangle |^2}{(1-\|x\|^2)^2}.
\ea

The authors also proved in \cite{BGM} that:
\ba \label{formulita}
\| \widetilde{\nabla} f(x) \|^2= (1-\|x\|^2)\left( \| \nabla f (x)\|^2 -|\R f(x)|^2 \right).
\ea

In order to simplify notation, for an analytic self-map $\psi: B_E \to B_E$, $x \in B_E$ and $w \in E$ we will denote:
\ba \label{ByC}
B(x,w)=\|\phi_{\psi(x)}'(0)^{-1}(\psi'(x)(w))\|
\ea
\ban
\mbox{and}  \ C(x,w)=\|\phi_x'(0)^{-1}(w)\|.
\ean

\noindent The following lemma just need easy calculations:

\begin{lema} \label{calculitos1}
Let $E$ be a complex Hilbert space, $\psi: B_E \to B_E$ an analytic self-map and $x \in B_E$. We have:
\begin{itemize} 
\item[a)] If $w \in E$ then:
\ba \label{ineq4}
\frac{\|w\|^2}{1-\|x\|^2} \leq C(x,w)^2 \leq \frac{\|w\|^2}{(1-\|x\|^2)^2}
\ea
and:
\ba \label{ineq5}
\frac{\|\psi'(x)(w)\|^2}{1-\|\psi(x)\|^2} \leq B(x,w)^2 \leq \frac{\|\psi'(x)(w)\|^2}{(1-\|\psi(x)\|^2)^2}
\ea

\item[b)] If there exists $w_x \in E$ satisfying $\psi'(x)(w_x)=\|\psi'(x)\| \psi(x)$ then:
\ba \label{ineq2}
\frac{\|\psi'(x)\|\|\psi(x)\|}{1-\|\psi(x)\|^2} = B(x,w_x) \leq \frac{\|\psi'(x)\|}{1-\|\psi(x)\|^2}
\ea	
and if, in addition, $w_x \neq 0$, then:
\ba \label{ineq3}
\ \ \ \ \ \ \ \frac{B(x,w_x)}{C(x,w_x)} \geq \tau_{\psi}(x) \frac{\|\psi(x)\|}{\|w_x\|}.
\ea
\end{itemize}
\end{lema}

\dem To prove a), by (\ref{ByC}) and (\ref{coci}) we have:
\ban
C(x,w)^2=\frac{(1-\|x\|^2) \|w\|^2+| \langle w,x\rangle|^2}{(1-\|x\|^2)^2}
\ean
\noindent so:
\ban
\frac{\|w\|^2}{(1-\|x\|^2)} \leq C(x,w)^2 \leq \frac{\|w\|^2}{(1-\|x\|^2)^2}
\ean

\noindent where last inequality is clear since $|\langle w,x \rangle| \leq \|w \| \|x\|$. Hence we conclude inequalities (\ref{ineq4}). The proof for (\ref{ineq5}) follows the same pattern. \smallskip

\noindent To prove b), making calculations and bearing in mind the expression of $B(x,w_x)$ in (\ref{ByC}) and (\ref{coci}) we have:
\ban
B(x,w_x)^2= \frac{(1-\|\psi(x)\|^2) \| \psi'(x)(w_x)\|^2+| \langle \psi'(x)(w_x),\psi(x)\rangle|^2}{(1-\| \psi(x)\|^2)^2} = \\ \frac{(1-\|\psi(x)\|^2)\|\psi'(x)\|^2  \| \psi(x)\|^2+\| \psi(x)\|^4 \|\psi'(x)\|^2}{(1-\| \psi(x)\|^2)^2}=\\ \frac{\|\psi'(x)\|^2  \| \psi(x)\|^2}{(1-\| \psi(x)\|^2)^2}
\ean
and we conclude inequality (\ref{ineq2}). This inequality together with inequality (\ref{ineq4}) result in inequality (\ref{ineq3}) since:
\ban
\frac{B(x,w_x)}{C(x,w_x)} \geq \frac{1-\|x\|^2}{1-\|\psi(x)\|^2} \frac{\|\psi'(x)\|\|\psi(x)\|}{\|w_x\|} 
\ean
and we are done. \qed \bigskip

From Lemma \ref{calculitos1} we easily conclude:
\begin{lema} \label{lema_ult}
For any $x \in B_E$ and $w \in E \setminus \{0\}$:
\ba \label{ineq6}
\frac{B(x,w)}{C(x,w)} \leq \frac{\sqrt{1-\|x\|^2}}{1-\|\psi(x)\|^2} \left\| \psi'(x) \left( \frac{w}{\|w\|}\right)\right\|
\ea
\noindent and:
\ba \label{ineq7}
\frac{B(x,w)}{C(x,w)} \geq \frac{1-\|x\|^2}{\sqrt{1-\|\psi(x)\|^2}} \left\| \psi'(x) \left( \frac{w}{\|w\|}\right)\right\|
\ea
\end{lema}

\medskip

We will also need the following result:
\begin{lema} \label{ineq_h}
Let $\psi: B_E \to B_E$ be an analytic map. Then for any $x \in B_E$ and $w \in E \setminus \{0\}$ we have:
$$\frac{B(x,w)}{C(x,w)}=\frac{\|\phi_{\psi(x)}'(0)^{-1}(\psi'(x)(w))\|}{\|\phi_{x}'(0)^{-1}(w)\|} \leq 1.$$
\end{lema}

\dem First suppose that $\psi(0)=0$ and consider the analytic function $f: \D \to \C$ given by:
$$f(z)=\frac{\langle \psi(z w) , \psi'(0)(w) \rangle}{\| \psi'(0)(w)\|}.$$
\noindent We suppose without loss of generality that $w$ belongs to the unit spehere $S_E$ of $E$ since otherwise we can divide both numerator and denominator by $\|w\|$. It is clear that $|f(z)| < 1$ and $f(0)=0$. By the Schwarz lemma we have that $|f'(0)|=\| \psi'(0)(w) \| \leq 1$ or equivalently, $\| \psi'(0)(w) \| \leq \|w\|$ for any $w \in E$. Consider $\mu= \phi_{\psi(x)} \circ \psi \circ \phi_x$. Notice that $\mu$ is a well-defined analytic self-map on $B_E$ and $\mu(0)=\phi_{\psi(x)}(\psi(x))=0$, so $\|  \mu'(0)(w)\| \leq \|w\|$. So for any $w \in E$ we have:
$$\|\phi_{\psi(x)}'(\psi(x)) \circ \psi'(x) \circ \phi_x'(0)(w)\| \leq \|w\|$$
\noindent and since $\phi_x'(0)^{-1}$ is a bijection on $B_E$, take $v \in E$ such that $w=\phi_x'(0)^{-1}(v)$ and we obtain for any $v \in E$:
$$\| \phi_{\psi(x)}'(\psi(x)) \circ \psi'(x) (v) \| \leq \|\phi_x'(0)^{-1}(v)\|$$
\noindent By Lemma \ref{derivadas} we know that $\phi_{\psi(x)}'(\psi(x))=\phi_{\psi(x)}'(0)^{-1}$ and we conclude the result. \qed \bigskip

The following corollary generalizes a result given by Kalaj (see \cite{Ka}) to the infinite dimensional setting :
\begin{cor} \label{cor_ineq}
Let $E$ be a complex Hilbert space and $\psi: B_E \to B_E$ an analytic map. Then for any $x \in B_E$:
$$\frac{1-\|x\|^2}{\sqrt{1-\|\psi(x)\|^2}} \| \psi'(x)\| \leq 1$$
\end{cor}

\dem It is sufficient to apply Lemma \ref{ineq_h} and inequality (\ref{ineq7}) in Lemma \ref{calculitos1}. \qed \bigskip

\begin{rem}
Kalaj proved that Corollary \ref{cor_ineq} is sharp by considering for any $t \in (0,\pi/2)$ the analytic self-map $\psi_t: B_2 \to B_2$ given by $\psi_t(z,w)=( z \sin t, \cos t)$ (see \cite{Ka}).
\end{rem}

\subsubsection{\textbf{Main results on bounded below composition operators on $\B(B_E)$}}

We will apply results on the automorphisms $\phi_x$ to the study of bounded below composition operators. First we provide a necessary condition by adapting the proof of Theorem 2 in \cite{DJO}:
\begin{teo}
Let $E$ be a complex Hilbert space and $\phi: B_E \to B_E$ be an analytic map. Suppose that $C_{\psi}: \B(B_E) \to \B(B_E)$ is bounded below. Then there exists $\eps >0$ and $0 < r < 1$ such that if $y \in B_E$ we have $\rho(\phi(x_y),y) \leq r$ for any $x_y \in B_E$ satisfying $\widetilde{\tau_\psi}(x_y) \geq \eps$.  
\end{teo}

\dem Suppose that $C_\psi$ is bounded below. Let $y \in B_E$ and consider the analytic function $f: B_E \to \C$ given by $f_y(x)=1/(1-\langle x,y \rangle).$ \noindent It is easy that:
$$f_y'(x)=\frac{y}{(1-\langle x,y \rangle)^2}$$
\noindent so:
\ban
\|f_y\|_\B=\sup_{x \in B_E} (1-\|x\|^2) \|f_y'(x)\|=\sup_{x \in B_E} (1-\|x\|^2)\frac{\|y\|}{|1-\langle x,y \rangle|^2} = \\ \sup_{x \in B_E} \|y\|\frac{1-\|\phi_y(x)\|^2}{1-\|y\|^2}=  \frac{\|y\|}{1-\|y\|^2}.
\ean

The analytic function $g_y: B_E \to \C$ given by $g_y(x)=\displaystyle f_y(x)/\|f_y\|_\B$ satisfies $\|g_y\|_\I \geq \|g_y\|_\B=1$. By Lemma \ref{lema_simp}, there exists $k >0$ such that $\|g_y \circ \psi \|_\I \geq k \|g_y\|_\I$ so bearing in mind that:
$$\| g_y \circ \psi\|_\I=\sup_{x \in B_E} \| \widetilde{\nabla} (g_y \circ \psi)(x)\|,$$ there exists $x_y \in B_E$ such that  $\| \widetilde{\nabla} (g_y \circ \psi)(x_y) \| \geq k/2$. So: 
\ban
\frac{k}{2} \leq \| \widetilde{\nabla} (g_y \circ \psi)(x_y) \| = \sup_{w \in E \setminus \{0\} } \| \widetilde{\nabla} (g_y \circ \psi)(x_y) (w) \| = \\  \sup_{w \in E \setminus \{0\} } \frac{|g_y'(\psi(x_y))(\psi'(x_y)(w))|}{\|\phi_{x_y}'(0)^{-1}(w)\|} =
\ean
\begin{align} \label{casi_ult}
\sup_{w \in E \setminus \{0\} } \frac{|g_y'(\psi(x_y))(\psi'(x_y)(w))|}{B(x_y,w)} 
\frac{B(x_y,w)}{C(x_y,w)}
\leq \|\widetilde{\nabla} g_y (\psi(x_y))\| \widetilde{\tau_{\psi}}(x_y)
\end{align}

\noindent where last inequality is clear by  (\ref{310})  and (\ref{ineq6}) in Lemma \ref{lema_ult}. By (\ref{formulita}) we have that:
\ban
\|\widetilde{\nabla} g_y (\psi(x_y))\|^2=(1-\|\psi(x_y)\|^2)( \| \nabla g_y (\psi(x_y))\|^2-|\R g_y (\psi(x_y))|^2)=\\ (1-\|\psi(x_y)\|^2) \frac{(1-\|y\|^2)^2}{\|y\|^2}\left(\frac{\|y\|^2}{|1-\langle \psi(x_y),y \rangle|^4}-\frac{|\langle \psi(x_y),y\rangle|^2}{|1-\langle \psi(x_y),y \rangle|^4}\right)=\\ (1-\|\psi(x_y)\|^2)(1-\|y\|^2)^2 \frac{1-\left| \Big \langle \psi(x_y),\frac{y}{\|y\|} \Big \rangle \right|^2}{|1-\langle \psi(x_y),y \rangle|^4}.
\ean
\noindent Notice that $|1-\langle a,b/\|b\| \rangle| \leq 2 |1-\langle a,b \rangle|$ for any $a,b \in B_E$ because:
\ban
|1-\langle a,b/\|b\| \rangle| \leq |1-\langle a,b \rangle|+|\langle a,b-b/\|b\| \rangle| \leq \\ |1-\langle a,b \rangle|+ 1-\|b\| \leq |1-\langle a,b \rangle|+1-|\langle a,b \rangle| = 2 |1-\langle a,b \rangle|.
\ean
\noindent Since:
$$1-\left| \Big \langle \psi(x_y),\frac{y}{\|y\|} \Big \rangle \right|^2 \leq  \left(1+\left| \Big \langle \psi(x_y),\frac{y}{\|y\|} \Big \rangle \right| \right) \left(1-\left| \Big \langle \psi(x_y),\frac{y}{\|y\|} \Big \rangle \right| \right)$$ we have:
\ban
\|\widetilde{\nabla} g (\psi(x_y))\|^2 \leq 4 (1-\|\psi(x_y)\|^2)(1-\|y\|^2) \frac{1}{|1-\langle \psi(x_y),y \rangle|^2} = \\ 4 (1-\| \phi_y (\psi(x_y))\|^2)= 4(1-\rho(y,\psi(x_y))^2).
\ean

\noindent Hence:
\ban
\frac{k}{2} \leq 2 (1-\rho(y,\psi(x_y))^2)^{1/2} \widetilde{\tau_{\psi}}(x_y)
\ean
\noindent which is satisfied if and only if:
\ban
\frac{k}{4} \leq (1-\rho(y,\psi(x_y))^2)^{1/2} \widetilde{\tau_{\psi}}(x_y) \leq  \widetilde{\tau_{\psi}}(x_y)
\ean
\noindent and we conclude that $\widetilde{\tau_{\psi}}(x_y) \geq \frac{k}{4}$. \medskip

From (\ref{casi_ult}) we have:
$$\frac{k}{2} \leq 2 (1-\rho(y,\psi(x_y))^2)^{1/2} \sup_{w \in E \setminus \{0\}}\frac{B(x_y,w)}{C(x_y,w)}$$
and by Lemma \ref{ineq_h} we conclude:
$$(1-\rho(y,\psi(x_y))^2)^{1/2} \geq k/4$$ \noindent which is satisfied if and only if:
$$\rho(y,\psi(x_y)) \leq \sqrt{1-k^2/16}.$$
\noindent Notice that if we choose another $x_y'$ satisfying $\widetilde{\tau_{\psi}}(x_y') \geq \frac{k}{4}$, then we have $\rho(y,\psi(x_y')) \leq \sqrt{1-k^2/16}$ by following the same argument from (\ref{casi_ult}). Take $r=\sqrt{1-k^2/16}$ and $\eps=k/4$ and we are done. \qed \bigskip


Now we provide a sufficient condition for a composition operator $C_\psi: \B(B_E) \to \B(B_E)$ to be bounded below. In order to extend the classical result given in Theorem \ref{teo5}, we will add a condition: we will need that $\psi(x_y)$ belongs to the range of $\psi'(x_y)$.

\begin{teo} \label{teo12}
Let $E$ be a complex Hilbert space and $\psi: B_E \to B_E$ be an analytic map. Suppose that there exist constants $0 < r < \frac{1}{15 A_0}$ and $\eps > 0$  such that for each $y \in \B_E$ there is a point $x_y \in B_E$ satisfying $\rho(\psi(x_y),y)<r$ and ${\tau_\psi}(x_y) > \eps$. Suppose in addition that $\psi(x_y) \in \psi'(x_y)(E)$. Then $C_{\psi}: \B(B_E) \to \B(B_E)$ is bounded below.
\end{teo}

\dem Let $f \in \B(B_E)$ such that $\|f\|_\I=1$. We will prove that there exists $k >0$ such that $\|f \circ \psi \|_\I \geq k$. By (\ref{seminormas}) we have that  $\|f\|_\R \geq \|f\|_\I /A_0$ so $\|f\|_\R \geq 1/A_0$. Take $y \in B_E$ such that $|Rf(y)|(1-\|y\|^2) \geq 14/(15 A_0)$. Then there exists $x_y \in B_E$ such that $\rho(y,\psi(x_y)) < r$ and $\tau_\psi(x_y) > \eps$. 
By (\ref{invar}) and (\ref{310}) and bearing in mind (\ref{ByC}), we have that for any $w \in E \setminus \{0\}$:
\ban
\|f \circ \psi\|_\I =\sup_{x \in B_E} \| \widetilde{\nabla} (f \circ \psi)(x)\| \geq \\ \frac{|(f \circ \psi)'(x_y)(w)|}{\|\phi_{x_y}'(0)^{-1}(w)\|}= 
\frac{|f'(\psi(x_y))(\psi'(x_y)(w))|}{B(x_y,w)} \frac{B(x_y,w)}{C(x_y,w)}.
\ean
\noindent Since $\psi(x_y) \in \psi'(x_y)(E)$, there exists $w_x \in E$ such that $\psi'(x_y)(w_x)=\| \psi'(x_y)\| \psi(x_y)$ so, in particular, the inequality above is true if we take $w_x$. By (\ref{ineq2}) in Lemma \ref{calculitos1} we have:
\ban
\frac{|f'(\psi(x_y))(\psi'(x_y)(w_x))|}{B(x_y,w_x)} =\frac{|f'(\psi(x_y))(\|\psi'(x_y)\| \psi(x_y))|}{B(x_y,w_x)}=\\ \frac{\|\psi'(x_y)\| |f'(\psi(x_y))( \psi(x_y))|(1-\|\psi(x_y)\|^2)}{\|\psi'(x_y)\|  \| \psi(x_y)\|} = \frac{|\R f(
\psi(x_y))| (1-\|\psi(x_y)\|^2)}{\| \psi(x_y)\|}
\ean
\noindent so:
\ban
\|f \circ \psi\|_\I \geq \R f(\psi(x_y)) (1-\|\psi(x_y)\|^2)\frac{B(x_y,w_x)}{C(x_y,w_x)}
\ean
\noindent and by (\ref{ineq3}) in Lemma \ref{calculitos1} we get:
\ban
\|f \circ \psi\|_\I \geq  \frac{|\R f(\psi(x_y))| (1-\|\psi(x_y)\|^2)}{\| \psi(x_y)\|} \frac{\|\psi(x_y)\| \tau_\psi(x_y)}{\|w_x\|}  \geq \\ |\R f(\psi(x_y))| (1-\|\psi(x_y)\|^2) \frac{\eps}{\|w_x\|}.
\ean
\noindent Using Corollary \ref{cor45}, we get:
\ban
| |\R f(\psi(x_y))| (1-\|\psi(x_y)\|^2)-|\R f(y)| (1-\|y\|^2)| \leq  14 \|f\|_\I \rho_E(\psi(x_y),y)
\ean
\noindent and since $\|f\|_\I=1$, we have:
\ban
\|f \circ \psi\|_\I \geq (|\R f(y)| (1-\|y\|^2)| - 14 \rho(\psi(x_y),y)) \eps \geq \\ \left(\frac{14}{15 A_0} - 14 r \right) \frac{\eps}{\|w_x\|} 
\ean
\noindent so we take:
$$k=\displaystyle 14\left(\frac{1}{15 A_0} -  r \right) \frac{\eps}{\|w_x\|} >0$$ \noindent and we conclude $\|C_\psi(f)\|_\I \geq k$. \qed \bigskip

Finally we check that for any $a \in B_E$, the automorphisms $\phi_a$ of $B_E$ satisfy conditions of Theorem \ref{teo12}. First, we provide the following proposition, which asserts that $\tau_{\phi_a}(x) \geq 1$ for any $x \in B_E$.

\begin{lema} \label{lemita2}
If $a \in B_E$ then $\tau_{\phi_a}(x) \geq 1$ for any $x \in B_E$.
\end{lema}

\dem By (\ref{pseudoh}) we have:
\ban
\frac{1-\|x\|^2}{1-\| \phi_a(x)\|^2}= \frac{|1-\langle x,a \rangle|^2}{1-\|a\|^2}
\ean
\noindent and since $\phi_a(x)=(P_a +s_a Q_a) \left( m_a(x)\right)$, then:
\ban
\phi_a'(x)=(P_a+s_a Q_a)'(m_a(x)) \circ m_a'(x)=(P_a+s_a Q_a)(m_a'(x))
\ean
\noindent so:
\ban
\| \phi_a'(x)\|^2= \|P_a(m_a'(x))\|^2+s_a^2 \|Q_a(m_a'(x))\|^2 \geq \|P_a(m_a'(x))\|^2.
\ean
It is clear that:
\ban
m_a'(x)(y)=\frac{-(1-\langle x,a \rangle)y+ \langle y,a \rangle (a-x)}{(1-\langle x, a \rangle)^2}
\ean
\noindent so:
\ban
\|\phi_a'(x)\| \geq \|P_a(m_a'(x))\| = \sup_{y \in \overline{B}_E} \|P_a(m_a'(x))(y))\| \geq \\ \left\|P_a \left(m_a'(x) \left( \frac{a}{\|a\|}\right) \right) \right\|=
\left\|P_a \left(
\frac{-(1-\langle x,a \rangle) \frac{a}{\|a\|}+ \langle \frac{a}{\|a\|},a \rangle (a-x)}{(1-\langle x, a \rangle)^2} \right) \right\|
\ean

\noindent and we deduce:
\ban
\tau_{\phi_a}(x) \geq \frac{|1-\langle x,a \rangle|^2}{1-\|a\|^2} \left\|P_a \left(
\frac{-(1-\langle x,a \rangle) \frac{a}{\|a\|}+ \langle \frac{a}{\|a\|},a \rangle (a-x)}{(1-\langle x, a \rangle)^2} \right) \right\| = \\ \frac{1}{1-\|a\|^2} \left\| P_a \left(
-(1-\langle x,a \rangle) \frac{a}{\|a\|}+ \langle \frac{a}{\|a\|},a \rangle (a-x) \right) \right\|= \\ \frac{1}{1-\|a\|^2} \left\| \left(
-(1-\langle x,a \rangle) \frac{a}{\|a\|}+ \|a\| a-  \frac{\langle x,a \rangle}{\|a\|^2} \|a\| a \right) \right\|= \\ \frac{1}{1-\|a\|^2} \left\| -\frac{1-\|a\|^2}{\|a\|} a \right\|=1
\ean
\noindent and we get $\tau_{\phi_a}(x) \geq 1$ as we wanted. \qed \bigskip


\begin{rem}
For any $a \in B_E$, the automorphisms $\phi_a$ of $B_E$ satisfy conditions of Theorem \ref{teo12} since by Proposition \ref{lemita2} we have:
$$\frac{1-\|x\|^2}{1-\|\phi_a(x)\|^2} \| \phi_a'(x)\| \geq 1$$
\noindent so we can take $\eps=1$, $r=0$ and for any $y \in B_E$ we choose $x_y= \phi_a(y)$. In addition, $\phi_a(x_y)=\phi_a(\phi_a(y))=y \in \phi_a'(x_y)(E)$ since $\phi_a'(x_y)$ is an invertible operator on $E$.
\end{rem}

%
%

\end{document}